\documentstyle[twoside,fleqn,espcrc2]{article}

\newcommand{\AmS}{{\protect\the\textfont2
  A\kern-.1667em\lower.5ex\hbox{M}\kern-.125emS}}

\title{Can we look at the quantisation rules as constraints?}
\author{Ennio Gozzi\address{Department of Theoretical Physics, 
University of Trieste\\
Strada Costiera 11, Miramare-Grignano 34014, Trieste\\
and INFN, Trieste, Italy.}}

\begin{document}

\begin{abstract}
In this paper we explore the idea of looking at the Dirac quantisation 
conditions as  $\hbar$-dependent constraints on the tangent bundle to
 phase-space. Starting from the path-integral version of
{\em classical} mechanics and using the natural Poisson brackets
structure present in the cotangent bundle to the tangent bundle of phase-
space, we handle the above constraints using the standard theory of Dirac for
constrained systems. The hope is to obtain, as total Hamiltonian, the Moyal
operator of time-evolution and as Dirac brackets the Moyal ones. Unfortunately
the program fails indicating that something is missing. We put forward at the
end some ideas for future work which may overcome this failure.
\end{abstract}

\maketitle

\section{INTRODUCTION}

Quantum mechanics (QM) is witnessing a revival of interest these days 
due to a variety of problems which range from the black-hole evaporation 
to the old issue of  measurement,  
to the problem of  the wave-function of the universe
and, last but not least, to the corrections which gravity would
induce on QM.  Somehow people are slowly  realizing that some of the  secrets of 
QM are still with us and partly unsolved.

We thought that it would help if we could look 
at QM in a different way. In this report
we will try to do that by looking\cite{gozzi1}\cite{gozzi2}\cite{gozzi3} at the 
quantisation conditions as constraints on the
tangent bundle to phase-space. This idea\cite{gozzi1} springs naturally 
from a 
path-integral formulation of classical mechanics (CM)\cite{gozzi4} 
we proposed in 1986-1989. 

The idea was further rekindled by some work
\cite{jar1} in which the authors tried to rederive
geometric quantisation\cite{wood1} using the tools of constrained and gauge
systems. Along this direction further work appeared later\cite{fradk1}
, work in which the authors obtained a global formulation of 
geometric quantisation via the theory of constraints.

Our work\cite{gozzi1}\cite{gozzi2}\cite{gozzi3} instead was not aimed at 
obtaining geometric quantisation or
at studying some   global issues in it, but was aimed at 
obtaining the Weyl-Wigner-Moyal formulation of QM\cite{moyal1}. The reason
for trying to obtain that particular formulation was  because that one 
was the  formulation of QM most closely related to our CM path-integral
\cite{gozzi4}. Both in Moyal QM and in CM the states are distributions 
in phase-space and the
quantum ones (Wigner functions), being\cite{berry1}  
both\footnote{By ~$S^{(1)}$~ we mean that the integral over all of phase-space
of the distributions
(not of the norm) is finite.} ~$S^{(1)}$~ and\footnote{At least for pure
density matrices or for mixed ones with a finite number of components.} 
~$L^{(2)}$ , are a subset of   properly
{\it enlarged} classical states which are only ~$S^{(1)}$. So if the quantum ones
are a subset of the {\it enlarged}-classical distributions,
this strongly indicates the presence of
a mechanism similar the one of gauge theory where the physical Hilbert space
appears as a subset of the full space. 
Working instead with geometric quantisation\cite{jar1}\cite{fradk1}
\cite{jar1} one goes 
from classical distributions in phase-space which are only ~$S^{(1)}$~ to 
wave-functions which are $L^{(2)}$~ but not necessarily~$S^{(1)}$~ . 
So the mechanism of reduction is less transparent.

Besides these aspects, our method\cite{gozzi1}\cite{gozzi2}\cite{gozzi3}
seems more natural in the sense that the strange 
auxiliary fields, which are needed\cite{jar1}\cite{fradk1}, appear
quite naturally in our formulation of CM\cite{gozzi4} springing out of our
path-integral as the basis of the tangent bundle
to phase-space. Moreover, differently than in refs.
\cite{jar1}\cite{fradk1}, we never have to 
quantize at the end, but we should pass naturally from the classical path-integral
\cite{gozzi4} to the quantum one of evolution of Wigner functions
\cite{gozzi5}. Anyhow  the works\cite{jar1}\cite{fradk1} which lead
 from CM to 
geometric quantisation are very interesting ones and they have achieved
remarkable results. Our work\cite{gozzi1}\cite{gozzi2}\cite{gozzi3} instead, 
even if more natural, it still afflicted with some problems as will be 
explained in what follows. Problems which we hope to overcome in the near 
future.

\section{QUANTISATION CONDITIONS AS CONSTRAINTS.}

The Dirac quantisation conditions are the usual rules of quantisation:

\begin{equation}
p  \longrightarrow -i\hbar\frac{\partial}{\partial q}
\end{equation}

In technical language this condition can be interpreted as an ~
$\hbar$-dependent constraint between elements ~$\phi\equiv(q,p)$~of
the phase-space ~$\cal{M}$~of the system and elements
~$(\frac{\partial}{\partial p},\frac{\partial}{\partial q})$~of the tangent
space ~$T_{\phi}\cal{M}$ in ~$\phi$~to ~$\cal{M}$.
In order to really interpret this as a constraint one should  work
in a space which unifies the phase-space~$\cal{M}$~with all the
tangent spaces~$T_{\phi}\cal{M}$~of ~$\cal{M}$. This space is
what is known as the tangent bundle to phase space and it is indicated
by ~$T\cal{M}$. A formulation of classical mechanics (CM)  in ~
$T\cal{M}$~has been given in ref.\cite{gozzi4}. We will skip the details
here and ask the reader to consult that paper. The final result
is that CM can be given a path-integral formulation whose generating
functional is:
\begin{equation}
$$Z_{cl}=\int {\cal D}\mu~exp~i\int dt \bigl\{{\widetilde{\cal L}}+
source~terms\bigr\}
\end{equation}
where ~${\cal D}\mu\equiv
{\cal D}\phi^{a}(t)~{\cal D}\lambda_{a}(t)~{\cal D}
c^{a}(t)~{\cal D}{\bar c}_{a}$. The variables ~$\phi^{a}$~ 
are the phase-space coordinates while the ~$\lambda_{a}$~
$c^{a}$~${\bar c}_{a}$~ are respectively commuting and anticommuting
auxiliary variables. The Lagrangian is~
${\widetilde{\cal L}}=\lambda_{a}{\dot{\phi}}^{a}+i{\bar c}_{a}{\dot c}
^{b}-{\widetilde{\cal H}}$~with the "Hamiltonian" given by
\begin{equation}
{\widetilde{\cal H}}=\lambda_{a}h^{a}+i{\bar c}_{a}\partial_{b}h^{a}
c^{b}
\end{equation}
and where ~$h^{a}$~ are the components of the Hamiltonian vector
field\cite{marsden}~$h^{a}(\phi)\equiv\omega^{ab}\partial_{b}H(\phi)$.
The various auxiliary variables introduced above have all a precise geometrical
meaning\cite{gozzi4}. For example the ~$c^{a}$~ are a basis for the forms\cite{marsden} 
, so~$c^{a}\equiv d\phi^{a}$. This means that the coordinates
~$(\phi^{a},c^{a})$~make up the cotangent bundle to phase-space, i.e~$T^{\ast}{\cal
M}$. We
refer the reader for details to ref.\cite{gozzi4}. 
To find the geometrical meaning of the ~$\lambda_{a}$~ and the
~${\bar c}_{a}$, we should first note\cite{gozzi4} that from the 
path-integral of eq.(2), one could define and calculate graded commutators
among the various variables. The result is the following
~$\bigl[\phi^{a},\lambda_{b}\bigr] =  i\delta^{a}_{b}$~and
$\bigl[c^{a},{\bar c}_{b}\bigr]=  \delta_{b}^{a}$. This implies that
we can realize ~$\lambda_{a}$~ and ~${\bar c}_{a}$~operatorially as
~$\lambda_{a}=-i\frac{\partial}{\partial\phi^{a}}\equiv
-i\partial_{a}~~~;~~~{\bar c}_{a}=\frac{\partial}{\partial c^{a}}$. 
Noticing this, we can say that ~$\lambda_{a}$~and ~
${\bar c}_{a}$~make a basis of the tangent space to the ~$(\phi^{a}, c^{a})$
space. So the 8n-coordinates ~$(\phi^{a},c^{a},\lambda_{a},{\bar c}_{a})$~
are a basis for the ~$T(T^{\ast}{\cal M})$~ which is the tangent bundle to the
cotangent bundle to phase-space. Inserting the operatorial realization of
~$\lambda_{a}$~and~${\bar c}_{a}$ into ~${\widetilde{\cal H}}$~ we would get
 what is known\cite{marsden} as the Lie-derivative of the Hamiltonian flow\cite{gozzi4}
. The first term of which is the Liouville operator
of classical mechanics:~${\hat L}=\frac{\partial H}{\partial p}\frac
{\partial}{\partial q}-\frac{\partial H}{\partial q}\frac{\partial}
{\partial p}$. The usual quantisation rule of eq.(1), if inserted in
the Liouville operator above, would not make sense because we would have to
say what is~$\frac{\partial}{\partial p}$~ at the {\it quantum}
level where ~$p$~ has become an operator. This means knowning what it
is the derivative of an operator. This belongs to the realm of non-commutative
geometry.We will adopt a less sofisticated strategy. As both CM\cite{gozzi4} and 
QM\cite{gozzi5} are now formulated via-path-integrals we will try to go
from one to the other via the Dirac\cite{regge} theory of constraints.
The reader may object that the constraints usually act on a classical
phase-space and not in an operatorial space. Well, in ref.\cite{gozzi4}
we proved that our path-integral, besides providing an operatorial
version of CM, is also naturally endowed with a classical Poisson brackets
structure just because the space~$T(T^{\ast}{\cal M})$~is isomorphic
to~$T^{\ast}(T{\cal M})$~ which is a cotangent bundle. We will call this
Poisson structure an extended Poisson structure
~$(epb)$. The~$\{\cdot,\cdot\}_{epb}$~among the basic fields are
$\bigl\{\phi^{a},\lambda_{b}\bigr\}_{epb}=\delta^{a}_{b}$, ~$\bigl\{c^{a},
{\bar c}_{b}\bigl\}_{epb}=-i\delta^{a}_{b}$~ and they reproduce the
classical equations of motion via the Hamiltonian of eq.(3)
: $\{\phi^a,H\}_{pb}=\bigl\{\phi^{a},{\widetilde{\cal H}}\bigr\}_{epb}$
where~the$\{\cdot,\cdot\}_{pb}$~are the old standard Poisson brackets
on ~${\cal M}$, which are~$\{\phi^{a},\phi^{b}\}_{pb}=\omega^{ab}$.
Basically this formulation of classical mechanics provides a manner (
if we disregard the anticommuting variables)
to generate the dynamics on the full~$T{\cal M}$. In this space now, remembering
the operatorial meaning of ~$\lambda$, the quantisation conditions can be 
written as the following constraints:
\begin{equation}
\Phi^{a}_{(0)}\equiv\theta(a-n)\bigl(\phi^{a}+\hbar\omega^{ab}\lambda_{b}\bigr)
\end{equation}
where the ~$\theta(a-n)$~is the standard step function which tells us that
the first ~$n$ indeces ~$a$, which are indicating
the q-variables, are put to zero.
Now that we have a CM in this enlarged space, we should treat the above 
constraint as it is usually done in constrained theory.

The first thing to do is to calculate
the {\it secondary} constraints which are those obtained by the evolution of
the primary ones~$\Phi^{a}_{(0)}$. The sub-index ~"$(0)$"~is to indicate that
it is primary. We will use the sub-indeces ~"$(1),(2),(3),...$"~to indicate
the {\it secondary} and {\it tertiary}, etc. constraints generated this way.
This procedure starts by using what is called\cite{regge} the total Hamiltonian
~${\widetilde{\cal H}}_{T}$~defined as
~${\widetilde{\cal H}}_{T}\equiv {\widetilde{\cal H}}+\sum_{a}u_{a}
\Phi^{a}_{(0)}$
~where ~$u_{a}$~are Lagrange multipliers. The secondary constraint~$\Phi^{a}
_{(1)}$ are then generated by imposing that the primary ones are left
invariant by the time evolution under ~${\widetilde{\cal H}}_{T}$.
In this case only  secondary ones are generated because the next step
starts determing the Lagrange multiplier. The secondary constraints turn out 
to be:~$\Phi^{a}_{(1)}=\theta(a-n)\bigl(\omega^{ab}\frac{\partial H}
{\partial\phi^{b}}-\hbar\omega^{ab}\lambda_{e}\omega^{ef}
\frac{\partial^{2}H}{\partial\phi^{b}\partial\phi^{f}}\bigr)$. 
With the same procedure it is possible to determine the Lagrange multipliers
and we refer the interested reader to ref.\cite{gozzi3}.
The constraints are second class and the next step is to build the associated
Dirac brackets\cite{regge}(which we will call extended Dirac brackets ($edb$)
because they are based on the extended Poisson brackets) between two 
observables ~$F$~ and ~$G$. These brackets are defined as
$\{F,G\}_{edb}\equiv\{F,G\}_{epb}-\{F,\Psi_{\alpha}\}_{epb}C^{\alpha\beta}
\{\Psi_{\beta},G\}_{epb}$~
where we have indicated collectively with ~$\Psi_{\alpha}$~the set of
primary ~$\Phi^{a}_{(0)}$~and secondary~$\Phi^{a}_{(1)}$~constraints
and~$C^{\alpha\beta}$~is the inverse of the matrix~$C_{\alpha\beta}$~
defined as~
$C_{\alpha\beta}\equiv\{\Psi_{\alpha},\Psi_{\beta}\}_{epb}$.
Applying all this to an harmonic oscillator~$H=1/2(p^{2}+q^{2})
$~we obtain\cite{gozzi3}, for the constrained evolution of an observable ~$F$,
the following: $\{F,{\widetilde{\cal H}}_{T}\}_{edb}= p\frac{\partial F}{\partial q}-
q\frac{\partial F}{\partial p}$. For an Hamiltonian with quartic potential:
$H={(1/2)}p^{2}+{(1/4)}q^{4}$~ the analog extended Dirac bracket evolution
would give:$\{F,{\widetilde{\cal H}}_{T}\}_{edb}={3/2}p\frac{\partial F}{\partial q}
-q^{3}\frac{\partial F}{\partial p}$.

Now that we have this formulation, which we could call an~
${\hbar}$-{\it constrained CM}~($\hbar$-CCM),
we want to compare it with real QM. Of course
the formulation of QM, with which we want to compare our~$\hbar$-CCM,
must be a formulation in phase-space (as our ~$\hbar$-CCM is) and
moreover it  must handle everything not with operators but with c-numbers as
our ~$\hbar$-CCM does. This formulation of QM exists and it was provided
by Weyl,Wigner and Moyal\cite{moyal1}. We will not review here all the 
Moyal formalism. Let us say that it is a procedure which replaces all
operators of standard QM with functions on phase-space (called symbols
of the operator) and the commutators
with brackets called Moyal brackets~$(mb)$. So the Heisenberg evolution of
the density matrix~$i\hbar~\partial_{t}{\widehat\varrho}=-\bigl[{\widehat\varrho},{\widehat
H}\bigr]$~ goes into~$\partial_{t}\varrho(\phi^{a},t)=-\bigl\{\varrho,H\bigr\}_{mb}$
where the~$\{\cdot,\cdot\}_{mb}$ are the Moyal brackets which are defined
as: 

$\bigl\{A,B\bigr\}_{mb}  =   A(\phi)~\frac{2}{\hbar}sin\bigl[
\frac{\hbar}{2}\stackrel{\leftarrow}{\partial_{a}}\omega^{ab}\stackrel
{\rightarrow}{\partial_{b}}\bigr]B(\phi)=\bigl\{A,B\bigr\}_{pb}+O(\hbar^{2})$.
In the classical limit ~($\hbar\rightarrow 0$)~the Moyal brackets reduces
to the classical Poisson brackets. The next thing to ask is if, 
as we did in CM, we can lift the action of
the Moyal Bracket on the ~$T{\cal M}$. The answer is yes\cite{gozzi5} and the result
looks as follows: It is possible to build an ~${\widetilde{\cal H}}^{\hbar}$
and a set of extended Moyal brackets~$\{\cdot,\cdot\}_{emb}$ such that:
~$\{\phi^{a},H\}_{mb}=\{\phi^{a},
{\widetilde{\cal H}}^{\hbar}\}_{emb}$. For the interested
reader the details are written in ref.\cite{gozzi5}. For the purpose of this
report it is enough to mention here the explicit form of
~${\widetilde{\cal H}}^{\hbar}$, which is:

${\widetilde{\cal H}}^{\hbar}={\widetilde{\cal H}}+{(\hbar)^{2}\over
3!}{\cal M}_{(1)}+{(\hbar)^{4}\over 5!}{\cal M}_{(2)}+\cdots$~where
\begin{equation}
{\cal M}_{(j)}  = \overbrace{[\lambda_{a}\omega^{ab}][\lambda_{c}\omega^{cd}]
\cdots[\lambda_{e}\omega^{ef}]}^{j+2}[\partial_{b}\partial_{d}\cdots
\partial_{f}H]
\end{equation}
The idea now is to compare the quantum evolution, provided by the ~
${\widetilde{\cal H}}^{\hbar}$~ under the extended Moyal brackets, with the
$\hbar$-CCM time evolution provided by the ~${\widetilde{\cal H}}_{T}$.
It is easy to see\cite{gozzi3} that the two evolution are the same for the 
harmonic oscillator, but they are different for the quartic potential. 
Even for the harmonic  oscillator the correspondence of QM with our $\hbar$-CCM is not perfect:
the extended Moyal bracket between two observables is not the same
as the extended Dirac bracket between the same observables unless one
of the two is  quadratic in~$\phi^{a}$. So our program of generating QM
from CM via a constraint has failed. Let us see which ones may be the reasons.
One which is immediately evident is the fact that we applied
the Dirac procedure to an Hamiltonian~${\widetilde{\cal H}}$~ which is not
singular. What we should have done is to start from an Hamiltonian
which already contains the quantisation constraints~~${\cal H}\equiv{\widetilde{\cal H}}
+\xi_{a}\Phi^{a}_{(0)}$~ and consider the Lagrange multipliers ~$\xi_{a}$
as dynamical variables. Anyhow even doing this  does not help
us going from~${\widetilde{\cal H}}$~ to~${\widetilde{\cal H}}^{\hbar}$. The reason
is that~${\widetilde{\cal H}}^{\hbar}$~ contains an infinite number of extra
pieces (see eq.(5)) beyond the first one which is~${\widetilde{\cal H}}$.
The idea put forward in ref.\cite{gozzi3} was that one should have generated
an infinite set of secondary and tertiary constraints. This of course does not happen.
What is more likely instead is that the constraints~$\Phi^{a}_{(0)}$~ are not
the full story. They may just be the first term of a constraint with
infinite terms (a formal power serie in $\hbar$). It this were so
, then it would be natural to expect that the total Hamiltonian
itself, ~${\widetilde{\cal H}}_{T}$, become a formal power serie in ~$\hbar$.
Having got to this point, the reader expert in constraints must have realized
that actually, even if we manage to build the ~${\widetilde{\cal H}}_{T}$,
its effect via the $(edb)$~is the same as the one of~${\widetilde{\cal H}}$:
$\{(\cdot),{\widetilde{\cal H}}_{T}\}_{edb}=\{(\cdot),{\widetilde{\cal H}}\}_{edb}$.
So what we have to check is if there are ~$(edb)$~such that
~$\{(\cdot),{\widetilde{\cal H}}\}_{edb}=\{(\cdot), 
{\widetilde{\cal H}}^{\hbar}\}_{emb}$.
Again the expert reader may immediately object that Dirac-Poisson
brackets can never be Moyal brackets for the different way they
behave with respect to derivation~\cite{lic78}. We feel hat theorem 
can actually be bypassed if the constraint is not a polinomial but a formal 
infinite series as in our case. The fact that it is such an infinite series
may puzzle the reader.  Actually it is easy to check that
eq.(1) is the correct quantisation rule if we work on wave-functions
and Hilbert space. It is not the correct one if we work with Wigner
functions which is the space of states we want to obtain. The Wigner 
function is defined as follows~$\varrho(p,q)\equiv\frac{1}{2\pi}\int
\psi^{\ast}(q-1/2\hbar s)e^{-ips}\psi(q+1/2\hbar s)ds$. Applying
on it the LHS of eq.(1) we get something different than applying
the RHS. What we get is: 
$\propto\int\bigl[\mp\frac{\partial\psi^{\ast}}{\partial q}\psi+\psi^{\ast}\frac{\partial
\psi}{\partial q}\bigr]e^{-isp}ds$~ where the ~$\mp$~ signs refer respectively
to the action of the LHS or the RHS of eq.(1).
This indicates that the constraint~$\Phi^{a}_{(0)}$ has to be modified.
One last thing before concluding. If we think that there is a constraint
operating here , then we must envision that the space of states 
(Wigner-functions)
be obtainable as a subset of the classical set of states. In order to do that,
we have to allow the classical states to be {\it non necessarily} positive
functions in phase-space and even discontinuous but still~$S^{(1)}$~ and to be
also distributions depending eventually even on external parameters.
It this enlargement (which is allowed by our CM path-integral) is done, then 
it is easy to prove that the Wigner functions, which are  both\cite{berry1} 
~$S^{(1)}$~~and ~$L^{(2)}$, are really a subset of the {\it enlarged} classical ones.

\end{document}